\documentclass[a4study,twoside,11pt]{article}
\usepackage{amssymb}
\usepackage{amsmath}
\usepackage{graphicx}
\usepackage{amsthm}
\usepackage{cite}
\usepackage{array}

\usepackage{amsmath,latexsym,amssymb,amsfonts,amsbsy,amsthm}
\usepackage{graphicx}
\usepackage{subfigure}
\usepackage[a4paper]{geometry}
\usepackage{epstopdf}
\usepackage{diagbox}
\usepackage{enumerate}
\usepackage{color}

\usepackage{comment}
\usepackage{indentfirst}

\newfont{\bb}{msbm10}

\newtheorem{theorem}{Theorem}[section]

\usepackage{verbatim}
\usepackage{algorithm}
\usepackage{algorithmicx}
\usepackage{algpseudocode}
\usepackage{amsmath}
\usepackage{graphics}
\usepackage{epsfig}
\usepackage{arydshln}

\usepackage{mathrsfs}
\usepackage{graphicx}
\usepackage{subfigure}
\usepackage{caption}

\usepackage{multirow} 
\usepackage{threeparttable}
\numberwithin{equation}{section}

\usepackage{float}

\usepackage{footnote}
\makesavenoteenv{table}

\usepackage{bm}
\usepackage{hyperref}
\hypersetup{
	hidelinks,
	colorlinks=true,
	linkcolor=blue,
	citecolor=blue
}

\title{On fast greedy block Kaczmarz methods for solving large consistent linear systems }
\author{
	A-Qin Xiao\\
	School of Mathematical Sciences, Tongji University,\\
	Shanghai, 200092, PR China.\\
	Email:xiaoaqin@tongji.edu.cn\\
	Jun-Feng Yin\thanks{Corresponding author.}\\
	School of Mathematical Sciences, Tongji University,\\
	Shanghai, 200092, PR China.\\
	Email:yinjf@tongji.edu.cn\\
	and\\
	Ning Zheng\\
	School of Mathematical Sciences, Tongji University,\\
	Shanghai, 200092, PR China.\\
	Email:nzheng@tongji.edu.cn\\
}
\date{}

\begin{document}
\cleardoublepage \pagestyle{myheadings}
	
\markboth{\small A.-Q. Xiao, J.-F. Yin and N. Zheng}
{\small On fast greedy block Kaczmarz methods for solving large consistent linear systems }
\captionsetup[figure]{labelfont={bf},labelformat={default},labelsep=period,name={Fig.}}	
\captionsetup[table]{labelfont={bf},labelformat={default},labelsep=period,name={Table }}		
\maketitle

\begin{abstract}
 A class of fast greedy block Kaczmarz methods combined with general greedy strategy and average technique are proposed for solving large consistent linear systems. Theoretical analysis of the convergence of the proposed method is given in detail. Numerical experiments show that the proposed methods are efficient and faster than the existing methods.\\


\end{abstract}

\section{Introduction}\label{sec1}
Consider the solution of consistent linear algebraic equations
\begin{equation}\label{eq1system}
Ax=b,
\end{equation}
where
$A \in \mathbb{R}^{m \times n}$ and $b \in \mathbb{R}^{m}$,
one of the classical and popular iteration methods is the Kaczmarz method \cite{1937K}. Due to its simplicity and efficiency, it was deeply studied and widely used in many practical scientific and engineering applications, for instance, computer tomography(CT) \cite{2001KS}, image reconstruction \cite{2008HD}, machine learning \cite{2016NSW} and option pricing \cite{2019FGNS}.

Let $A^{(i)}$ be the $i$-th row of matrix $A$ and $b^{(i)}$ be the $i$-th entry of vector $b$, respectively.
Given an initial guess vector $x_0\in\mathbb{R}^n$, the classical Kaczmarz method iterates by
\begin{equation*}
x_{k+1}= x_k +\frac{b^{(i_k)}-A^{(i_k )}x_k }{\left\| A^{(i_k )}\right\|^2_2 }(A^{(i_k )})^T,  \quad k= 0,1,2,\cdots,
\end{equation*}
where $i_k=(k \text{ mod } m)+1$. 
To improve the convergence of Kaczmarz method, Strohmer and Vershynin \cite{2009SV} proposed a randomized Kaczmarz method by selecting the row index $i_k$ with probability proportion to $\left\| A^{(i_k )}\right\|^2_2$ and proved its linear convergence rate in expectation.
Bai and Wu \cite{2018BW} constructed a greedy randomized Kaczmarz method to
accelerate the convergence performance. For more variants of the randomized and greedy Kaczmarz methods, we refer the reader to \cite{2018BWrgrk,2019BW,2021DHZMY,2022YLZ}.

The idea of block Kaczmarz method can date back to the work of Elfving in\cite{1980E}, which used many equations simultaneously at each iteration. The block Kaczmarz method can be described as
\begin{equation}\label{eqc1BK}
x_{k+1}= x_k + A_{\mathcal{J}_k}^{\dagger}(b_{\mathcal{J}_k}- A_{\mathcal{J}_k}x_k),\quad k=0,1,2,\cdots,
\end{equation}
where $A_{\mathcal{J}_k}^{\dagger}$ represents the Moore-Penrose pseudoinverse of the chosen submatrix $A_{\mathcal{J}_k}$ and $\mathcal{J}_k$ is the block row indices. Needell and Tropp \cite{2014NT} proposed a randomized block Kaczmarz method, where $\mathcal{J}_k$  is chosen uniformly at random from $[m]=\left\lbrace 1,2,\dots,m \right\rbrace $.
Further, Niu and Zheng \cite{2020NZ} proposed a greedy block Kaczmarz method which adaptively choose the block row indices without predetermining a partition of the row indices of the matrix $A$.

However, each iterate of the  block Kaczmarz method require the computation of the pseudoinverse of the selected submatrix corresponding to the residual subvector and it usually costs expensive. Necoara \cite{2019N} established a unified framework of randomized average block Kaczmarz methods by taking a convex combination of several updatings as a new direction and implemented on the distributed computing units. Miao and Wu \cite{2022MW} proposed an average block variant of the greedy randomized Kaczmarz method \cite{2018BW} and studied its convergence. 

In this work, we construct a class of fast greedy block Kaczmarz methods with average technique to avoid computing the pseudoinverse of submatrices of the coefficients matrix, where a modified greedy strategy utilizing the general norm of residual vectors is proposed and well studied, which can choose the working rows based on this greedy criterion dynamically and flexibly.
Theoretical analysis of the convergence of the proposed method is given in detail.
Numerical experiments further demonstrate that the proposed methods are efficient and faster than the existing methods.

The rest of this paper is organized as follows. In Section \ref{sec2},
a class of fast greedy block Kaczmarz methods is presented with
a modified greedy row selection strategy and the convergence theory of the proposed method is established.
Numerical experiments are reported in Section \ref{sec3}  to display the efficiency of the new method. Finally, we draw the conclusions  in Section \ref{sec4}.

\section{The fast greedy block Kaczmarz method}\label{sec2}
In this section, after reviewing the fast deterministic block Kaczmarz method,
we present a class of fast greedy block Kaczmarz methods for solving large consistent linear systems \eqref{eq1system} by using a modified greedy row selection strategy and the averaging technique.

Let $\xi_k$ be a linear combination of unit column vectors $e_i\in \mathbb{R}^m$ $(i\in\tau_k)$ and its coefficients are the corresponding entries of the residual vector, that is,
\begin{equation*}
 \xi_k=\sum\limits_{i \in \mathcal{I}_{k}}\left(b^{(i)}-A^{(i)} x_{k}\right) e_{i},\quad k= 0,1,2,\cdots.
\end{equation*}
Similar to the iteration of average block Kaczmarz methods \cite{2019N},
the stepsize  and weight are set to be $\alpha_k= \frac{\|\xi_k\|^2_2 \| A_{\mathcal{I}_{k}} \|^2_F }{\| A^T \xi_k \|^2_2}$ and  $\omega_i^k= \frac{\| A^{(i)}\|_2^2}{\| A_{\mathcal{I}_{k}}\|_F^2}$, $i\in \mathcal{I}_{k}$ respectively, the fast deterministic block Kaczmarz method can iterate as follows
 \begin{align}
 \notag
 x_{k+1}&= x_k + \alpha_k \left(\sum\limits_{i\in \mathcal{I}_{k}}\omega_i^k\frac{ b^{(i)}-A^{(i)}x_k }{ \| A^{(i)}\|_2^2 } (A^{(i)})^T \right)\\ \notag
 &= x_k +  \frac{\|\xi_k\|^2_2 \| A_{\mathcal{I}_{k}} \|^2_F }{\| A^T \xi_k \|^2_2}\left(\sum\limits_{i\in \mathcal{I}_{k}}\frac{\| A^{(i)}\|_2^2}{\| A_{\mathcal{I}_{k}}\|_F^2}\cdot\frac{ b^{(i)}-A^{(i)}x_k }{ \| A^{(i)}\|_2^2 } (A^{(i)})^T \right)\\ \notag
 &= x_k +  \frac{\xi_k^T(b-Ax_k) }{\| A^T \xi_k \|^2_2}\left(\sum\limits_{i\in \mathcal{I}_{k}}\left(  b^{(i)}-A^{(i)}x_k\right) (A^{(i)})^T \right)\\ \label{eqfdbkinduc}
 &=x_{k}+\frac{\xi_{k}^{T}\left(b-A x_{k}\right)}{\left\|A^{T} \xi_{k}\right\|_{2}^2 } A^{T} \xi_{k}.
 \end{align}
where the block indices $\mathcal{I}_{k}$ is chosen by
$$
\mathcal{I}_{k} = \left\{i \big| \lvert b^{(i)}-A^{(i)} x_{k}\rvert^2 \geq \gamma_{k}
\left\|b- Ax_k\right\|_2^2\left\|A^{(i)}\right\|_2^2 \right\}
$$
with
$$\gamma_{k}= \frac{1}{2}\left(\frac{1}{\left\|b- Ax_k\right\|_2^2} \max\limits_{1 \leq i \leq m}\left\{\frac{\lvert b^{(i)}-A^{(i)} x_{k}\rvert^{2}}{\left\|A^{(i)}\right\|_2^2}\right\}+ \frac{1}{\left\|A\right\|_F^2} \right). $$

One drawback of the fast deterministic block Kaczmarz method is the greedy strategy for choosing the working block indices. The size of the block $\mathcal{I}_{k}$ may be small if the Frobenius norm of the coefficient matrix $A$ is very small, which may lead to a very slow convergence.

In order to further accelerate  the fast deterministic block Kaczamarz method, we propose to determine the working rows by a modified greedy strategy that utilizes the general norm of residual vectors. Let $\eta \in(0,1]$ and $p\in[1,+\infty)$,  the control index subset is defined as
$$
\tau_{k}= \left\{i \big| \lvert b^{(i)}-A^{(i)} x_{k}\rvert^p \geq \epsilon_{k}\left\|A^{(i)}\right\|_p^p \right\}, \quad k=0,1,2\ldots,
$$
where
$$
\epsilon_{k}= \eta\cdot\max\limits_{1 \leq i \leq m}\left\{\frac{\lvert b^{(i)}-A^{(i)} x_{k}\rvert^{p}}{\left\|A^{(i)}\right\|_p^p}\right\}. 
$$
It follows that
\begin{equation*}
\epsilon_{k}= \eta\cdot \max _{1 \leq i \leq m}\left\{\frac{\lvert b^{(i)}-A^{(i)} x_{k}\rvert^{p}}{ \left\|A^{(i)}\right\|_p^p }\right\}\leq  \max _{1 \leq i \leq m}\left\{\frac{\lvert b^{(i)}-A^{(i)} x_{k}\rvert^{p}}{\left\|A^{(i)}\right\|_p^p }\right\},
\end{equation*}
which implies that there is at least one index $j\in[m] $ such that
\begin{equation*}
\frac{\lvert b^{(j)}-A^{(j)} x_{k}\rvert^{p}}{\left\|A^{(j)}\right\|_p^p}= \max _{1 \leq i \leq m}\left\{\frac{\lvert b^{(i)}-A^{(i)} x_{k}\rvert^{p}}{\left\|A^{(i)}\right\|_p^p}\right\},
\end{equation*}
then $j\in\tau_k$, i.e. $\tau_k$ is not empty.

Given an initial guess vector $x_0\in\mathbb{R}^n$, the fast greedy block Kaczmarz method is described in Algorithm \ref{algFGBK}.
\begin{algorithm}[!htbp]
	\caption{The fast greedy block Kaczmarz method (FGBK)} \label{algFGBK}
	\begin{algorithmic}[1]
		\Require $A, b, x_0,l$, $\eta\in (0,1]$ and $p\in[1,+\infty)$
		\Ensure  $ x_{l}$
		\For  {$k=0,1,2,\ldots,l-1$ }
		\State Compute
		\begin{equation*}
		\epsilon_{k}= \eta\cdot\max_{1 \leq i \leq m}\left\{\frac{\lvert b^{(i)}-A^{(i)} x_{k}\rvert^{p}}{\left\|A^{(i)}\right\|_p^p }\right\}.
		\end{equation*}
		\State  Determine the control index set of positive integers
		\begin{equation}\label{eqfgbk2}
		\tau_{k}= \left\{i\mid \lvert b^{(i)}-A^{(i)} x_{k}\rvert^p \geq \epsilon_{k}\left\|A^{(i)}\right\|_p^p \right\}.
		\end{equation}	
		\State Compute
		\begin{equation}\label{eqfgbk3}
		\xi_{k}=\sum\limits_{i \in \tau_{k}}\left(b^{(i)}-A^{(i)} x_{k}\right) e_{i}.
		\end{equation}
		\State Set
		\begin{equation}\label{eqfgbk4}
		x_{k+1}=x_{k}+\frac{\xi_{k}^{T}\left(b-A x_{k}\right)}{\left\|A^{T} \xi_{k}\right\|_{2}^2 } A^{T} \xi_{k}.
		\end{equation}
		\EndFor
	\end{algorithmic}
\end{algorithm}

Moreover, the theoretical analyses for the convergence performance of the fast greedy block Kaczmarz method are established as follows.

\begin{theorem}\label{themconvergence}
 Given an initial vector $x_0\in \mathbb{R}^{n}$ in the column space of $A^T$.
  Then, the  sequence $\left\{x_k \right\}_{k=0}^{\infty}$ generated by the fast greedy block Kaczmarz method converges to the unique least-norm solution $x_{\ast}= A^{\dagger}b$. Moreover, the norm of the error of the approximate solution satisfies
\begin{equation*}
\left\|x_{k+1}-x_{*}\right\|_{2}^{2} \leq  \left( 1-\beta_{k}(\eta,p)   \sigma^2_{\min }\left(A\right) \right)  \left\|x_{k}-x_{*}\right\|_{2}^{2}, \quad k \geq 0,
\end{equation*}
where
$\beta_{k}(\eta,p)= \frac{\eta^{\frac{2}{p}}}{\sum\limits_{i\in[m]\backslash\tau_{k-1}} \left\| A^{(i)}\right\|_p^{2}  } \cdot \frac{\sum\limits_{i \in \tau_{k}}  \left\|A^{(i)} \right\|^{2}_p }{\sigma^2_{\max }\left(  A_{\tau_{k}}\right)},$  $\eta\in (0,1]$ and $p\in[1,+\infty)$.
\begin{proof}
\rm For each $k>0$, denote the projector $P_{k}=\frac{A^{T} \xi_{k} \xi_{k}^{T} A}{\left\|A^{T} \xi_{k}\right\|_{2}^{2}}$,
since $P^T_k=P_k,P^2_k=P_k$, so $P_k$ is an orthogonal projection.
From the iterate scheme \eqref{eqfgbk4}, it holds that
\begin{equation*}
\begin{aligned}
x_{k+1}-x_{*} &=x_{k}-x_{*}+\frac{\xi_{k}^{T}\left(b-A x_{k}\right)}{\left\|A^{T} \xi_{k}\right\|_{2}^{2}} A^{T} \xi_{k} \\
&=x_{k}-x_{*}-\frac{\xi_{k}^{T} A\left(x_{k}-x_{*}\right)}{\left\|A^{T} \xi_{k}\right\|_{2}^{2}} A^{T} \xi_{k} \\
&=x_{k}-x_{*}-\frac{A^{T} \xi_{k} \xi_{k}^{T} A}{\left\|A^{T} \xi_{k}\right\|_{2}^{2}}\left(x_{k}-x_{*}\right).
\end{aligned}
\end{equation*}
By the Pythagorean theorem, it follows that
\begin{equation}\label{eqproof35}
\begin{aligned}
\left\|x_{k+1}-x_{*}\right\|_{2}^{2} &=\left\|\left(I-P_{k}\right)\left(x_{k}-x_{*}\right)\right\|_{2}^{2  } \\
&=\left\|x_{k}-x_{*}\right\|_{2}^{2}-\left\|P_{k}\left(x_{k}-x_{*}\right)\right\|_{2}^{2} \\
&=\left\|x_{k}-x_{*}\right\|_{2}^{2}-\left\|\frac{\xi_{k}^{T} A\left(x_{k}-x_{*}\right)}{\left\|A^{T} \xi_{k}\right\|_{2}^{2}} A^{T} \xi_{k}\right\|_{2}^{2} \\
&=\left\|x_{k}-x_{*}\right\|_{2}^{2}-\frac{\lvert \xi_{k}^{T}\left(b-A x_{k}\right)\rvert  ^{2}}{\left\|A^{T} \xi_{k}\right\|_{2}^{2}}.
\end{aligned}
\end{equation}
Let $E_{k} \in \mathbb{R}^{m \times\lvert \tau_{k}\rvert  }$ be the matrix whose columns  are consisted of all the vector $e_{i} \in \mathbb{R}^{m}$ with $i \in \tau_{k}$. Denote $A_{\tau_{k}}=E_{k}^{T} A$, $\widehat{\xi}_{k}=E_{k}^{T} \xi_{k}$, then
\begin{equation}\label{eqproof36}
\left\|\widehat{\xi}_{k}\right\|_{2}^{2}=\xi_{k}^{T} E_{k} E_{k}^{T} \xi_{k}=\left\|\xi_{k}\right\|_{2}^{2}=\sum_{i \in \tau_{k}}\lvert b^{(i)}-A^{(i)} x_{k}\rvert  ^{2},
\end{equation}
and
\begin{equation}\label{eqproof37}
\left\|A^{T} \xi_{k}\right\|_{2}^{2}=\xi_{k}^{T} A A^{T} \xi_{k}=\widehat{\xi}_{k}^{T} E_{k}^{T} A A^{T} E_{k} \widehat{\xi}_{k}=\widehat{\xi}_{k}^{T} A_{\tau_{k}} A_{\tau_{k}}^{T} \widehat{\xi}_{k}=\left\|A_{\tau_{k}}^{T} \widehat{\xi}_{k}\right\|_{2}^{2}.
\end{equation}
Therefore,
\begin{equation}\label{eqproof38}
\left\|A_{\tau_{k}}^{T} \widehat{\xi}_{k}\right\|_{2}^{2}=\widehat{\xi}_{k}^{T} A_{\tau_{k}} A_{\tau_{k}}^{T} \widehat{\xi}_{k} \leq \sigma^2_{\max }\left(A_{\tau_{k}}\right)\left\|\widehat{\xi}_{k}\right\|_{2}^{2},
\end{equation}
where $\sigma_{\max}(A_{\tau_k})$ represents the largest singular value of selected submatrix $A_{\tau_k}$.
By the definition of $\xi_{k}$ in \eqref{eqfgbk3} and \eqref{eqproof36}, it holds that

\begin{equation}\label{eqproof39}
\begin{aligned}
\xi_{k}^{T}\left(b-A x_{k}\right) &=\left(\sum_{i \in \tau_{k}}\left(b^{(i)}-A^{(i)} x_{k}\right) e_{i}^{T}\right)\left(b-A x_{k}\right) \\
&=\sum_{i \in \tau_{k}}\left(\left(b^{(i)}-A^{(i)} x_{k}\right) e_{i}^{T}\left(b-A x_{k}\right)\right) \\
&=\sum_{i \in \tau_{k}}\lvert b^{(i)}-A^{(i)} x_{k}\rvert  ^{2} \\
&=\left\|\widehat{\xi}_{k}\right\|_{2}^{2} .
\end{aligned}
\end{equation}
Since both $x_{k}$ and $x_{*} \in R(A^T)$,  $x_{k}-x_{*} \in \mathcal{R}\left(A^{T}\right)$, then
\begin{equation}\label{eqproof310}
\left\|b-A x_{k}\right\|_{2}^{2}=\left\|A\left(x_{k}-x_{*}\right)\right\|_{2}^{2} \geq \sigma^2_{\min }\left(A \right)\left\|x_{k}-x_{*}\right\|_{2}^{2}.
\end{equation}
From the equalities \eqref{eqproof37}--\eqref{eqproof310} and the definition of $\tau_k$ in \eqref{eqfgbk2}, it follows that
\begin{equation}\label{eqproof311}
\begin{aligned}
\frac{\lvert \xi_{k}^{T}\left(b-A x_{k}\right)\rvert  ^{2}}{\left\|A^{T} \xi_{k}\right\|_{2}^{2}} &=\frac{\left(\sum\limits_{i \in \tau_{k}}\lvert b^{(i)}-A^{(i)} x_{k}\rvert  ^{2}\right)\left\|\widehat{\xi}_{k}\right\|_{2}^{2}}{\left\|A_{\tau_{k}}^{T} \widehat{\xi}_{k}\right\|_{2}^{2}} \\
& \geq \frac{\sum\limits_{i \in \tau_{k}}\lvert b^{(i)}-A^{(i)} x_{k}\rvert  ^{2}}{\sigma^2_{\max }\left(  A_{\tau_{k}}\right)} \\
& \geq \frac{\sum\limits_{i \in \tau_{k}}\left(  \lvert b^{(i)}-A^{(i)} x_{k}\rvert  ^{p}\right)^{\frac{2}{p}}}{\sigma^2_{\max }\left(  A_{\tau_{k}}\right)} \\
&\geq(\epsilon_{k})^{\frac{2}{p}}\cdot\frac{ \sum\limits_{i \in \tau_{k}}\left(\left\|A^{(i)}\right\|^{p}_p\right)^{\frac{2}{p}}}{\sigma^2_{\max }\left(  A_{\tau_{k}}\right)} \\
&\geq(\epsilon_{k})^{\frac{2}{p}}\cdot\frac{ \sum\limits_{i \in \tau_{k}} \left\|A^{(i)}\right\|^{2}_p }{\sigma^2_{\max }\left(  A_{\tau_{k}}\right)} .
\end{aligned}
\end{equation}
In addition, it is seen that
\begin{align*}
b- Ax_{k}
&= b- A\left( x_{k-1}+ \frac{\xi_{k-1}^T(b_{\tau_{k-1}}-A_{\tau_{k-1}}x_{k-1}) }{\left\| A_{\tau_{k-1}}^T \xi_{k-1} \right\|^2_2 }A_{\tau_{k-1}}^T \xi_{k-1}\right) \\
&=\left(b-Ax_{k-1}\right)-\frac{A\xi_{k-1}^T(b_{\tau_{k-1}}-A_{\tau_{k-1}}x_{k-1}) }{\left\| A_{\tau_{k-1}}^T \xi_{k-1} \right\|^2_2 }A_{\tau_{k-1}}^T \xi_{k-1},k=1,2,\cdots.
\end{align*}
Therefore,
\begin{align*}
b_{\tau_{k-1}}- A_{\tau_{k-1}}x_{k} =& \left( b_{\tau_{k-1}}-A_{\tau_{k-1}}x_{k-1}\right)\\
&-\frac{A_{\tau_{k-1}}\xi_{k-1}^TA_{\tau_{k-1}}^T \xi_{k-1}}{\left\| A_{\tau_{k-1}}^T \xi_{k-1} \right\|^2_2} (b_{\tau_{k-1}}-A_{\tau_{k-1}}x_{k-1})= 0.
\end{align*}
It is known that
\begin{equation*}
\begin{aligned}
\|b-Ax_k\|_{2}^{2}&= \sum_{i\in[m]\backslash \tau_{k-1}} \frac{\lvert b^{(i)}-A^{(i)} x_{k}\rvert  ^{2}}{\left\|A^{(i)}\right\|_p^{2}} \left\|A^{(i)}\right\|_p^{2}\\
&= \sum_{i\in[m]\backslash \tau_{k-1}}\left(\frac{\lvert b^{(i)}-A^{(i)} x_{k}\rvert  ^{p}}{\left\|A^{(i)}\right\|_p^{p}  }\right)^{\frac{2}{p}} \left\|A^{(i)}\right\|_p^{2}\\
&\leq \left(\max _{1 \leq i \leq m}\left\{\frac{\lvert b^{(i)}-A^{(i)} x_{k}\rvert  ^{p}}{\left\|A^{(i)}\right\|_p^{p}}\right\}\right)^{\frac{2}{p}}\cdot\sum\limits_{i\in[m]\backslash\tau_{k-1}} \left\| A^{(i)}\right\|_p^{2}.
\end{aligned}
\end{equation*}
It follows that
\begin{equation}\label{eqproof312}
\begin{aligned}
(\epsilon_{k})^{\frac{2}{p}}&=\left(\eta\cdot \max _{1 \leq i \leq m}\left\{\frac{\lvert b^{(i)}-A^{(i)} x_{k}\rvert^{p}}{\left\|A^{(i)}\right\|_p^{p}}\right\}\right)^{\frac{2}{p}}
\geq \eta^{\frac{2}{p}} \cdot \frac{\|b-Ax_k\|_{2}^{2}} {\sum\limits_{i\in[m]\backslash\tau_{k-1}} \left\| A^{(i)}\right\|_p^{2} }  \\
&\geq \frac{\eta ^{\frac{2}{p}}}{\sum\limits_{i\in[m]\backslash \tau_{k-1}} \left\|A^{(i)}\right\|_p^{2} } \cdot \sigma^2_{\min } (A) \left\|x_{k}-x_{*}\right\|_{2}^{2}.
\end{aligned}
\end{equation}
From \eqref{eqproof311} and \eqref {eqproof312}, it deduces that
\begin{equation}\label{eqpfend}
\frac{\lvert \xi_{k}^{T}\left(b-A x_{k}\right)\rvert  ^{2}}{\left\|A^{T} \xi_{k}\right\|_{2}^{2}} \geq \beta_{k}(\eta,p)\cdot \sigma^2_{\min } (A)\left\|x_{k}-x_{*}\right\|_{2}^{2}.
\end{equation}
where $\beta_{k}(\eta,p)= \frac{\eta^{\frac{2}{p}}}{\sum\limits_{i\in[m]\backslash\tau_{k-1}} \left\| A^{(i)}\right\|_p^{2}  } \cdot \frac{\sum\limits_{i \in \tau_{k}}  \left\|A^{(i)} \right\|^{2}_p }{\sigma^2_{\max }\left(  A_{\tau_{k}}\right)}$, $\eta\in(0,1]$ and $p\in[1,+\infty)$.

\noindent Finally, by combining \eqref{eqproof35} and \eqref{eqpfend}, it follows that
\begin{equation*}
\left\|x_{k+1}-x_{*}\right\|_{2}^{2}\leq \left(1- \beta_{k}(\eta,p) \sigma^2_{\min } (A)\right)\left\|x_{k}-x_{*}\right\|_{2}^{2}.
\end{equation*}
\end{proof}
\end{theorem}
 Note that the upper bound of convergence rate of the fast greedy block Kaczmarz method is related to the relaxation parameter $\eta$,  the parameter $p$, the geometric properties of the coefficient matrix $A$ and its row submatrices at each iteration.  
 However, the practical convergence speed of the fast greedy block Kaczmarz methods could be faster than the upper bound.

\section{Numerical experiments} \label{sec3}
In this section, a number of numerical experiments are presented to illustrate the efficiency of the fast greedy block Kaczmarz (FGBK) method, compared with the greedy block Kaczmarz (GBK) method \cite{2020NZ} and the fast deterministic block Kaczmarz (FDBK) method \cite{2022CH} in terms of the number of iteration steps (denoted as `IT') and the elapsed computing time in seconds (denoted as `CPU').

In the numerical experiment, the solution vector $ x$ is firstly constructed and $b= Ax$ so that the linear system is consistent.
All the iterations are started from  the initial vector $x_0=0$, and terminated when the relative solution error (denoted as `RSE') satisfies
\begin{equation*}
\text{RSE}=\frac{\lVert x_k -x_{\ast} \rVert_2^2}{\lVert x_0- x_{\ast} \rVert_2^2}<10^{-6},
\end{equation*}
or the number of iteration steps exceeds a maximal number, e.g., 10000.
For the greedy block Kaczmarz method, the control index set is determined by 
$$
\mathcal{J}_{k}= \left\{i \Big| \lvert b^{(i)}-A^{(i)} x_{k}\rvert^2 \geq \delta_{k}\max\limits_{1 \leq i \leq m}\left\{\frac{\lvert b^{(i)}-A^{(i)} x_{k}\rvert^{2}}{\left\|A^{(i)}\right\|_2^2}\right\}\left\|A^{(i)}\right\|_2^2 \right\}
$$
with parameter $\delta_k$ is
$$
\delta_k= \frac{1}{2} + \frac{1}{2}\frac{\lVert b-Ax_k \rVert_2^2}{\lVert A \rVert_F^2}\left( \max_{1 \leq i \leq m}\left\{\frac{\lvert b^{(i)}-A^{(i)} x_{k}\rvert^{2}}{\left\|A^{(i)}\right\|_2^2 }\right\} \right)^{-1}
$$
while the relaxation parameter $\eta_{exp}$ in the fast greedy block Kaczmarz method is experimentally selected by minimizing the numbers of total iterations.

In the first example,  the matrix $A$ is  generated by the MATLAB function `randn',
where the size of the  matrices is chosen to be
$5000\times10000, 5000\times12000, 5000\times14000, 5000\times16000$ and $5000\times18000$, respectively.
In Table \ref{tab:resultrand}, the number of iterations and elapsed CPU time of the greedy block Kaczmarz, the fast deterministic block Kaczmarz and fast greedy block Kaczmarz methods with $p=1,2$ and $3$ are listed respectively.

From Table \ref{tab:resultrand}, it can observed that the GBK, FDBK, FGBK($p=1$), FGBK($p=2$) and FGBK($p=3$) methods converge successfully and the FGBK-type methods outperform the other two methods in terms of both the iteration count and CPU time. 
Moreover, the fast greedy block Kaczmarz method with $p=1$ requires the least number of iterations.
It shows that the efficiency of the modified greedy row selection strategy and indicates that a small value of $p$ may further accelerate the convergence.

\begin{table}[!htbp]
	\centering 
	\caption{ Numerical results for random matrices. }\label{tab:resultrand}
	\begin{tabular}{lllllll}
		\hline{Method}&{$ m \times n $} & 5000$\times$10000 & 5000$\times$12000 & 5000$\times$14000 & 5000$\times$16000 & 5000$\times$18000 \\
		\hline 
		\multirow{3}{*} {GBK}    & IT&	 543&	 349&	 251&	 208&	 169 \\  
		& CPU&	 40.3710&	 36.0279&	 33.5738&	 35.6595&	 35.6040 \\  
		& RSE&	 9.97$\times 10^{-7}$&	 9.92$\times 10^{-7}$&	 9.89$\times 10^{-7}$&	 9.82$\times 10^{-7}$&	 9.76$\times 10^{-7}$ \\  
		\hline
		\multirow{3}{*} {FDBK}      & IT&	 559&	 356&	 256&	 209&	 170 \\  
		& CPU&	 27.8735&	 21.4744&	 17.9981&	 16.5299&	 15.4575 \\  
		& RSE&	 9.91$\times 10^{-7}$&	 9.85$\times 10^{-7}$&	 9.71$\times 10^{-7}$&	 9.54$\times 10^{-7}$&	 9.97$\times 10^{-7}$ \\  
		\hline
		\multirow{4}{*} {FGBK(p=1)} & $\eta_{exp}$&	 0.10&	 0.05&	 0.05&	 0.05&	 0.05 \\  
		& IT&	 73&	 47&	 35&	 29&	 24 \\  
		& CPU&	 4.5057&	 3.6111&	 3.1248&	 2.9545&	 2.7597 \\  
		& RSE&	 9.53$\times 10^{-7}$&	 9.52$\times 10^{-7}$&	 9.29$\times 10^{-7}$&	 9.70$\times 10^{-7}$&	 8.80$\times 10^{-7}$ \\  
		\hline  
		\multirow{4}{*} {FGBK(p=2)} & $\eta_{exp}$&	 0.05&	 0.05&	 0.05&	 0.05&	 0.05 \\  
		& IT&	 74&	 48&	 36&	 30&	 25 \\  
		& CPU&	 4.0051&	 3.1378&	 2.8106&	 2.6807&	 2.4948 \\  
		& RSE&	 8.77$\times 10^{-7}$&	 8.03$\times 10^{-7}$&	 8.02$\times 10^{-7}$&	 7.92$\times 10^{-7}$&	 7.84$\times 10^{-7}$ \\  
		\hline  
		\multirow{4}{*} {FGBK(p=3)} & $\eta_{exp}$&	 0.05&	 0.05&	 0.05&	 0.05&	 0.05 \\  
		& IT&	 82&	 55&	 42&	 35&	 30 \\  
		& CPU&	 6.0167&	 5.6553&	 5.5993&	 6.0705&	 6.3334 \\  
		& RSE&	 9.36$\times 10^{-7}$&	 9.38$\times 10^{-7}$&	 8.38$\times 10^{-7}$&	 8.77$\times 10^{-7}$&	 9.37$\times 10^{-7}$ \\  
		\hline
	\end{tabular}
\end{table}

In the second exmaple, the matrices are taken from the SuiteSparse Matrix Collection \cite{2011DH}
to further compare the convergence performances of these Kaczmarz methods. The test matrices `stat96v5' and `crew1' come from linear programming problems while `bibd\_17\_8' and `bibd\_16\_8' come from combinatorial problems.
In Table \ref{tab:informarket}, the sizes ($m\times n$), rank, density and condition number of the test matrices are listed respectively.

\begin{table}[!htbp]
	\centering
	\caption{Information of the matrices from SuiteSparse Matrix Collection.}\label{tab:informarket}
	\begin{tabular}{lllll}
		\hline {Name} &stat96v5 &crew1 &bibd\_17\_8 &bibd\_16\_8\\ \hline	
		$ m \times n $ &2307 $\times$ 75779  &135 $\times$ 6469   &136 $\times$ 24310 &120 $\times$ 12870   \\
		rank &2307 &135  &136 &120  \\
		density &0.13\%  &5.38\%  &20.59\% &23.33\%   \\
		cond($A$) &19.52 &18.20 &9.04  &9.54    \\ \hline	
	\end{tabular}
\end{table}

In Table \ref{tab:resultmarket}, the number of iterations and CPU time of the greedy block Kaczmarz method, the fast deterministic block Kaczmarz method and the fast greedy block Kaczmarz method with $p=1,2$ and 3 are reported respectively.

 \begin{table}[!htbp]
	\centering
	\caption{ Numerical results for the matrices from SuiteSparse Matrix Collection.}\label{tab:resultmarket}
		\begin{tabular}{lllllll}
		\hline\multicolumn{2}{l} {Method} &stat96v5 &crew1  &bibd\_17\_8 & bibd\_16\_8 \\
		\hline 
		\multirow{3}{*} {GBK}    & IT&	 80&	 547&	 237&	 280 \\  
		& CPU&	 1.2582&	 0.7068&	 19.7138&	 11.3608 \\  
		& RSE&	 8.69$\times 10^{-7}$&	 9.85$\times 10^{-7}$&	 9.75$\times 10^{-7}$&	 9.83$\times 10^{-7}$ \\  
		\hline
		\multirow{3}{*} {FDBK}      & IT&	 249&	 815&	 256&	 289 \\  
		& CPU&	 0.2940&	 0.1438&	 0.3103&	 0.1846 \\  
		& RSE&	 9.22$\times 10^{-7}$&	 9.89$\times 10^{-7}$&	 9.66$\times 10^{-7}$&	 9.91$\times 10^{-7}$ \\  
		\hline
		\multirow{4}{*} {FGBK(p=1)} & $\eta_{exp}$&	 0.05&	 0.20&	 0.10&	 0.10 \\  
		& IT&	 36&	 356&	 125&	 138 \\  
		& CPU&	 0.1787&	 0.1065&	 0.1829&	 0.1249 \\  
		& RSE&	 8.46$\times 10^{-7}$&	 9.57$\times 10^{-7}$&	 9.98$\times 10^{-7}$&	 8.98$\times 10^{-7}$ \\  
		\hline  
		\multirow{4}{*} {FGBK(p=2)} & $\eta_{exp}$&	 0.05&	 0.30&	 0.15&	 0.15 \\  
		& IT&	 39&	 406&	 137&	 163 \\  
		& CPU&	 0.1859&	 0.0907&	 0.1735&	 0.1057 \\  
		& RSE&	 8.27$\times 10^{-7}$&	 9.93$\times 10^{-7}$&	 9.44$\times 10^{-7}$&	 9.33$\times 10^{-7}$ \\  
		\hline  
		\multirow{4}{*} {FGBK(p=3)} & $\eta_{exp}$&	 0.05&	 0.15&	 0.05&	 0.05 \\  
		& IT&	 45&	 387&	 134&	 163 \\  
		& CPU&	 0.2124&	 0.0945&	 0.1802&	 0.1137 \\  
		& RSE&	 9.26$\times 10^{-7}$&	 9.53$\times 10^{-7}$&	 9.17$\times 10^{-7}$&	 8.02$\times 10^{-7}$ \\  
		\hline  
	\end{tabular}  
\end{table}

From Table \ref{tab:resultmarket}, it is observed that the proposed methods require fewer iterations and less CPU time than the other two Kaczmarz methods.
For the matrix `stat96v5', the proposed method with $p=1$ has the least number of iteration and the least CPU time. For the other three matrices, the proposed method with $p=1$ requires the least number of iterations. It further indicates that a small value of $p$ may further improve the speed of convergence.

\begin{figure}[!htbp]
	\centering
	\subfigure[stat96v5]
	{
		\begin{minipage}[t]{0.4\linewidth}
			\centering
			\includegraphics[width=1\textwidth]{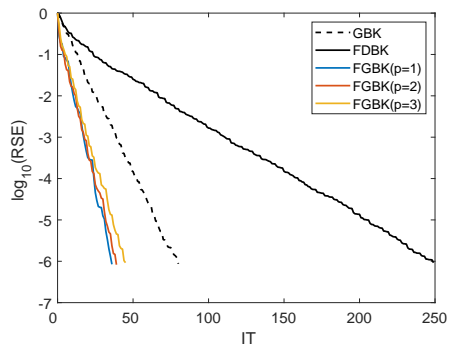}
		\end{minipage}
	}
	\subfigure[crew1]
	{
		\begin{minipage}[t]{0.4\linewidth}
			\centering
			\includegraphics[width=1\textwidth]{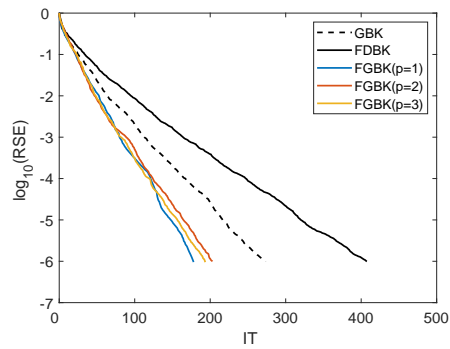}
		\end{minipage}
	}	
	\centering
	\subfigure[bibd\_17\_8]
	{
		\begin{minipage}[t]{0.4\linewidth}
			\centering
			\includegraphics[width=1\textwidth]{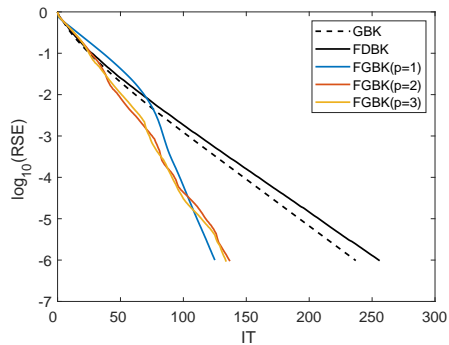}
		\end{minipage}
	}
	\subfigure[bibd\_16\_8]
	{
		\begin{minipage}[t]{0.4\linewidth}
			\centering
			\includegraphics[width=1\textwidth]{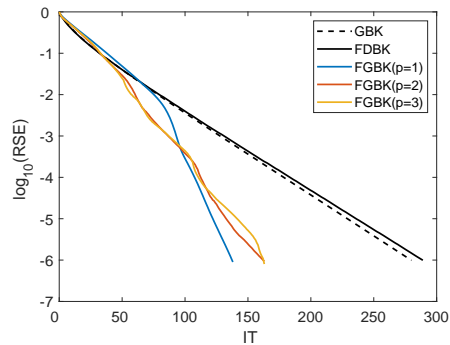}
		\end{minipage}
	}
	\caption{Convergence curves for the matrices from SuiteSparse Matrix Collection. }
	\label{fig:RSEvsITmarket}
\end{figure}

In Figure \ref{fig:RSEvsITmarket}, the curves of the relative solution error versus the number of iterations are plotted for GBK, FDBK, FGBK($p=1$), FGBK($p=2$) and FGBK($p=3$) respectively.

From Figure \ref{fig:RSEvsITmarket}, it is obviously seen that the fast greedy block Kaczmarz methods converge faster than the greedy block Kaczmarz method and the fast deterministic block Kaczmarz method, which confirms the numerical result in Table \ref{tab:resultmarket} and further shows the efficiency of the modified greedy row selection strategy.

\section{Conclusions}\label{sec4}
A class of fast greedy block Kaczmarz methods is presented for solving large consistent linear systems. Theoretical analysis proves the convergence of the proposed methods and show that the upper bound of the convergence rate is related to the geometric properties of the coefficient matrix and its block submatrices. Numerical experiments further illustrate that the proposed methods are efficient and faster than the fast deterministic block Kaczmarz method.

\vspace{1em}
\noindent{\bf Acknowledgements}  This work is supported by the National Natural Science Foundation of China (Grant No. 11971354).

\bibliographystyle{unsrt}

\end{document}